\documentstyle[11pt]{article}
\begin{document}
\author{{Shiqiu Zheng$^{1, 2}$\thanks{Corresponding author, E-mail: shiqiumath@163.com(S. Zheng).}\ , \ \ Shoumei Li$^{1}$\thanks{E-mail: lisma@bjut.edu.cn(S.Li).}}
  \\
\small(1, College of Applied Sciences, Beijing University of Technology, Beijing 100124, China)\\
\small(2, College of Sciences, Hebei United University,  Tangshan 063009, China)\\
}
\date{}

\title{\textbf{Representation theorems for generators of BSDEs with monotonic and convex growth generators}\thanks{This work is supported by the Science and Technology Program of Tangshan (No. 13130203z).}}\maketitle

\textbf{Abstract:}\quad In this paper, we establish representation theorems for generators of backward
stochastic differential equations (BSDEs in short), whose generators are monotonic and convex growth in $y$ and quadratic growth in $z$. We also obtain a converse comparison theorem for such BSDEs.
\\



\section{Introduction}
Representation theorem for generator of BSDE shows that the generator of BSDE can be represented by the limit of solutions of corresponding BSDEs. It is firstly established by Briand et al. (2000) for BSDEs whose generators satisfy Lipschitz condition and two additional assumptions that $E[\sup_{0\leq t\leq T}|g(t,0,0)|^2] <\infty$ and $(g(t, y, z))_{t\in[0,T]}$ is continuous
in $t$. Then it is generalized step by step by a series of work of Jiang (see Jiang (2005a, b, c, 2006, 2008))  to case that $g$ only satisfies Lipschitz condition. Since then, representation theorem for generators of BSDEs is further studied for more general generators by many papers. For example, Liu et al. (2007) studies it under so called Mao's non-Lipschitz condition. Jia (2008) and Fan $\&$ Jiang (2010)  studies it under continuous and linear growth condition, respectively. Ma $\&$ Yao (2010)  studies it under quadratic growth condition in $z$ and two additional assumptions. Fan $\&$ Jiang (2011)  studies it for generators which are monotonic and polynomial growth in $y$ and linear growth in $z$. Recently, Zheng (2014a) generalizes the representation theorem in Ma $\&$ Yao (2010) to the general case that generators only are linear growth in $y$ and quadratic growth in $z$. The main result of this paper is that we establish representation theorems for generators of BSDEs whose generators are monotonic and convex growth in $y$ and quadratic growth $z$.

The motivation of this paper is to understand the relation between generators and solutions of BSDEs considered in Briand $\&$ Hu (2008), which shows that BSDE has at least a solution when the generator $g$ is monotonic in $y$ and quadratic growth in $z.$  For this purpose, we want to establish a representation theorem for generators of BSDEs considered in Briand $\&$ Hu (2008). In this paper, we further make a assumption that generator $g$ is convex growth in $y$, which generalizes the polynomial growth condition considered in Fan $\&$ Jiang (2011). This convex growth condition can help us obtain Lemma 3.3,  which can be used to deal with the difficulty arising from the proof of representation theorem for generators of such BSDEs.

It is worth noting that representation theorems for generators is a powerful tool to interpret the relation between generators and solutions of BSDEs. It has played an important role in studying the properties of generators by virtue of solutions of corresponding BSDEs and nonlinear expectation theory. One can see Briand et al. (2000), Jiang (2005a, b, c, 2006, 2008), Jia (2008), Fan $\&$ Jiang (2010),  Ma $\&$ Yao (2010), Fan et al. (2011) and Zheng (2014a), etc.

This paper is organized as follows. In Section 2, we will give some basic notions. In Section 3, we will establish representation theorems for generators and a converse comparison theorem of BSDEs whose generators are monotonic and convex growth in $y$ and quadratic growth $z$.
\section{Preliminaries}
Let $(\Omega ,\cal{F},\mathit{P})$ be a complete probability space
carrying a $d$-dimensional standard Brownian motion ${{(B_t)}_{t\geq
0}}$, starting from $B_0=0$, let $({\cal{F}}_t)_{t\geq 0}$ denote
the natural filtration generated by ${{(B_t)}_{t\geq 0}}$, augmented
by the $\mathit{P}$-null sets of ${\cal{F}}$, let $|z|$ denote its
Euclidean norm, for $\mathit{z}\in {\mathbf{R}}^n$, let $T>0$ be a given real number. We define the following usual
spaces:

${\cal{L}}^1([0,T])=\{f(t):$ Lebesgue measurable function such that $\int_0^T|f(t)|dt<\infty\};$

$L^p({\mathcal {F}}_T)=\{\xi:\ {\cal{F}}_T$-measurable
random variable; $\|\xi\|_{L^p}=\left({\mathbf{E}}\left[|\xi|^p\right]\right)^{\{\frac{1}{p}\wedge1\}}<\infty\},\ \ \ p>0; $

$L^\infty({\mathcal {F}}_T)=\{\xi:{\cal{F}}_T$-measurable
random variable; $\|\xi\|_\infty=\textrm{esssup}_{\omega\in\Omega}|\xi|<\infty\};$

${\mathcal{S}}^p_T
 \left(\mathbf{R}\right)=\{\psi:$ continuous predictable
process; $\|\psi\|_{{\cal{S}}^{{p}}}^p={\mathbf{E}} \left
[{\mathrm{sup}}_{0\leq t\leq T} |\psi _t|^p\right] <\infty \},\ \ \ p\geq1;$

${\mathcal{S}}^\infty
 _T({\mathbf{R}})=\{\psi:$ continuous predictable
process; $\|\psi\|_\infty=\textrm{esssup}_{(\omega,t)\in\Omega\times [0,T]}|\psi _t| <\infty\};$

${\cal{H}}^p_T({\mathbf{R}}^d)=\{\psi:$ predictable
process;
$\|\psi\|_{{\cal{H}}^p}=\left({\mathbf{E}}\left[\int_0^T|\psi_t|^pdt\right]\right)^{\{\frac{1}{p}\wedge1\}}
<\infty \},\ \ \ \ p>0.$\\

For any progressively measurable stochastic process $\{h_t\}_{t\in[0,T]},$ we maybe use the following assumption (B) in this paper.

\textbf{Assumption (B)}. There exist $\Omega'\subset\Omega$ such that $P(\Omega')=1$ and $\{h^2_t(\omega)\}_{\omega\in\Omega'}\subset {\cal{L}}^1([0,T])$ is uniformly integrable. \\

\textbf{Remark 2.1} Clearly, if $\{h_t\}_{t\in[0,T]}$ satisfies assumption (B), then we have $\|\int_0^Th_t^2 dt\|_\infty<\infty$ and $\lim\limits_{\varepsilon\rightarrow0^+}\|\{\int_t^{t+\varepsilon}h^2_tdt\}\|_\infty=0.$\\

Let us consider a function ${g}\left( \omega ,t,y,z\right) : \Omega \times [0,T]\times \mathbf{%
R\times R}^{\mathit{d}}\longmapsto \mathbf{R},$ such that
$\left(g(t,y,z)\right)_{t\in [0,T]}$ is progressively measurable for
each $(y,z)\in\mathbf{%
R\times R}^{\mathit{d}}$. In this paper, we consider the following BSDE (Pardoux $\&$ Peng (1990)):
$$Y_t=\xi +\int_t^Tg(s,Y_s,Z_s)
ds-\int_t^TZ_s\cdot dB_s,\ \ \ t\in[0,T],$$
usually called BSDE with parameter $(g,T,\xi),$ where $g$ is called generator, $\xi$ and $T$ are called terminal variable and terminal time, respectively.\\

Now, we introduce a stochastic differential equation (SDE). Suppose $b(\cdot,\cdot,\cdot):\Omega\times [0,T]\times \textbf{R}^m\mapsto \textbf{R}^m$ and $\sigma(\cdot,\cdot,\cdot):\Omega\times [0,T]\times \textbf{R}^m\mapsto \textbf{R}^{m\times d}$ satisfy the following two conditions:

(H1) there exists a constant $\mu\geq0$ such that $P$-$a.s.$, $ \forall t\in[0,T],\
\forall x,y\in \textbf{R}^{\mathit{m}},$
$$|{b}(t,x)-b(t,y) |+|\sigma(t,x)-\sigma(t,y) |\leq \mu|x-y|.$$

(H2) there exists a constant $\nu\geq0$ such that $P$-$a.s.$, $ \forall t\in[0,T],\
\forall x\in \textbf{R}^{\mathit{m}},$

$$|b(t,x)|+|\sigma(t,x)|\leq \nu\left(1+|x|\right).$$
Given $(t,x)\in[0,T[\times \textbf{R}^m$, by SDE theory, the following SDE:
$$\left\{
  \begin{array}{ll}
     X_s^{t,x}=x +\int_t^sb(u,X_u^{t,x})du+\int_t^s\sigma(u,X_u^{t,x})dB_u,\ \ \  s\in]t,T],\\
     X_s^{t,x}=x,\ \ \   s\in[0,t],
  \end{array}
\right.$$
has a unique $s$-continuous adapted solution $X_s^{t,x}.$\\

\textbf{Lemma 2.1} (Lebesgue Lemma, see Hewitt and Stromberg (1978, Lemma 18.4)) Let $f$ be a Lebesgue integrable function on the interval $[0,T]$. Then for almost every $t\in[0,T[$, we have
$$\lim\limits_{\varepsilon\rightarrow0^+}\frac{1}{\varepsilon}
\int_{t}^{t+\varepsilon}|f(u)-f(t)|ds=0.$$

%%%%%%%%%%%%%%%%%%%%%%%%%%%%%%%%%%%%%%%%%%%%%%%%%%%%%%%%%%%%%%%%
\section{Representation theorem for generators}
%%%%%%%%%%%%%%%%%%%%%%%%%%%%%%%%%%%%%%%%%%%%%%%%%%%%%%%%%%%%%%%%
In this section, we will establish representation theorems for generators of BSDEs whose generators are monotonic and convex growth in $y$ and quadratic growth in $z$. \\

\textbf{Assumption (A)}. There exist constants
$\beta\geq0, \gamma>0$, a progressively measurable non-negative stochastic process $\{\alpha_t\}_{t\in[0,T]}$ and a strictly increasing convex function $\varphi:\ {\textbf{R}}^+\longrightarrow {\textbf{R}}^+$ with $\varphi(0)=\varphi(0+)=0,$ such that $P-a.s.,$

(i) for each $t\in [0,T],\ (y,z)\longmapsto g (t, y, z)$ is continuous;

(ii) (monotonicity in $y$) for each $(t,z)\in[0,T]\times{\mathbf{R}}^d,$
$$\forall y\in{\mathbf{R}},\ \ y(g(t,y,z)-g(t,0,z))\leq\beta|y|^2;$$

(iii) (convex growth in $y$) for each $(t,y,z)\in [0,T]\times{\mathbf{ R\times R}}^{\mathit{d}},$
$$|{g}\left( t,y,z\right)|\leq \alpha_t+\varphi(|y|)+\frac{\gamma}{2}|z|^2.$$

By Briand and Hu (2008, Lemma 2), if $g$ assumptions (A) and $\|\int_0^T\alpha_tdt\|_\infty<\infty$, then for each $\xi\in L^\infty({\mathcal {F}}_T),$ BSDE with parameter $(g,T,\xi)$ has a maximal solution $(\overline{{Y}}_t,\overline{{Z}}_t)\in{\mathcal{S}}^\infty_T
({\mathbf{R}})\times{\cal{H}}^{2}_T({\mathbf{R}}^d)$ and a minimal solution $(\underline{{Y}}_t,\underline{{Z}}_t)\in{\mathcal{S}}^\infty_T
({\mathbf{R}})\times{\cal{H}}^{2}_T({\mathbf{R}}^d)$  in the sense that for any solution $({Y}_t,{Z}_t)$ of BSDE with parameter $(g,T,\xi)$, we have $\underline{{Y}}_t\leq{Y}_t\leq\overline{{Y}}_t.$ By Briand and Hu (2008, Proposition 1), we also have the following facts:
$$\|{Y}_t\|_\infty\leq
    e^{\beta T}\left(\|\xi\|_\infty+\left\|\int_0^T\alpha_tdt\right\|_\infty\right),\eqno(1)$$
and for each $\varepsilon\in]0,T-t]$, if $({Y}_{s}^{t+\varepsilon},{Z}_{s}^{t+\varepsilon})$ is a solution of BSDE with parameter $(g,t+\varepsilon,0),$ then
$$\sup_{t\leq s\leq t+\varepsilon}|Y_{s}^{t+\varepsilon}|\leq\frac{1}{\gamma}
\log \left(E\left[\exp\left(\gamma e^{\beta\varepsilon}\int_t^{t+\varepsilon}\alpha_rdr\right)|{\cal{F}}_t\right]\right)\leq e^{\beta\varepsilon}\left\|\int_t^{t+\varepsilon}\alpha_rdr\right\|_\infty.
\eqno(2)$$
Moveover, if $\{\alpha_t\}_{t\in[0,T]}$ satisfies assumption (B), then by (2) and Remark 2.1, we have
$$\lim\limits_{\varepsilon\rightarrow0^+}\left\|\sup_{t\leq s\leq t+\varepsilon}|Y_{s}^{t+\varepsilon}|\right\|_\infty=0.
\eqno(3)$$
and if $\|\alpha_t\|_\infty<\infty,$ then by (2), we have
$$\sup_{t\leq s\leq t+\varepsilon}|Y_{s}^{t+\varepsilon}|\leq \varepsilon e^{\beta\varepsilon}\|\alpha_t\|_\infty.\eqno(4)$$

\textbf{Lemma 3.1} Let $g$ satisfies assumption (A) and $\{\alpha_t\}_{t\in[0,T]}$ in (A) satisfies assumption (B), for each $t\in[0,T[$ and stopping time $\tau\in ]0,T-t],$ we have
$$\lim\limits_{\varepsilon\rightarrow0^+}\frac{1}{\varepsilon}E\left[\int_{t}^{t+\varepsilon\wedge\tau}
|{Z}_r^{t+\varepsilon\wedge\tau}|^2dr|{\cal{F}}_t\right]=0\ \ \
\textrm{and}\ \ \ \lim\limits_{\varepsilon\rightarrow0^+}\frac{1}{\varepsilon}E\left[\int_{t}^{t+\varepsilon\wedge\tau}
|{Z}_r^{t+\varepsilon\wedge\tau}|^2dr\right]=0,
$$
where $(Y_{s}^{t+\varepsilon\wedge\tau},Z_{s}^{t+\varepsilon\wedge\tau})$ is an arbitrary solution of BSDE with parameter $(g,t+\varepsilon\wedge\tau,0).$\\

\emph{Proof.} For each $t\in[0,T[,\ \varepsilon\in]0,T-t]$ and any stopping time $\tau\in]0,T-t],$ we consider the following BSDEs with parameter $(g,t+\varepsilon\wedge\tau,0)$
$${Y}_s^{t+\varepsilon\wedge\tau}=\int_{s}^{{t+\varepsilon\wedge\tau}}g(r,{Y}_r^{t+\varepsilon\wedge\tau},
{Z}_r^{t+\varepsilon\wedge\tau})dr-\int_{s}^{t+\varepsilon\wedge\tau}{Z}_r^{t+\varepsilon\wedge\tau}dB_r.$$
By (ii) in assumption (A), we have
$$yg(t,y,z)\leq yg(t,0,z)+\beta|y|^2\leq y(\alpha_t+\frac{\gamma}{2}|z|^2)+\beta|y|^2.$$
Applying It\^{o} formula to $|{Y}_s^{t+\varepsilon\wedge\tau}|^2$ for $s\in[t,t+\varepsilon\wedge\tau],$ and in view of above inequality, we can deduce
\begin{eqnarray*}
&&|{Y}_t^{t+\varepsilon\wedge\tau}|^2+\int_{t}^{t+\varepsilon\wedge\tau}|{Z}_r^{t+\varepsilon\wedge\tau}|^2dr\\
&=&2\int_{t}^{t+\varepsilon\wedge\tau}{Y}_r^{t+\varepsilon\wedge\tau}g(r,{Y}_r^{t+\varepsilon\wedge\tau},
{Z}_r^{t+\varepsilon\wedge\tau})dr-2\int_{t}^{t+\varepsilon\wedge\tau}{Y}_r^{t+\varepsilon\wedge\tau}{Z}_r^{t+\varepsilon\wedge\tau}dB_r\\
&\leq&2\int_{t}^{t+\varepsilon\wedge\tau}\alpha_r|{Y}_r^{t+\varepsilon\wedge\tau}|dr
+2\beta\int_{t}^{t+\varepsilon\wedge\tau}|{Y}_r^{t+\varepsilon\wedge\tau}|^2dr
+\gamma\int_{t}^{t+\varepsilon\wedge\tau}|{Y}_r^{t+\varepsilon\wedge\tau}||{Z}_r^{t+\varepsilon\wedge\tau}|^2dr\\
&&-2\int_{t}^{t+\varepsilon\wedge\tau}{Y}_r^{t+\varepsilon\wedge\tau}{Z}_r^{t+\varepsilon\wedge\tau}dB_r\ \ \ \ \ \ \ \ \ \ \ \ \ \ \ \ \ \ \ \ \ \ \ \ \ \ \ \ \ \ \ \ \ \ \ \ \ \ \ \  \ \ \ \ \ \ \ \ \ \ \ \ \ \ \  \ \ \ \ \ \ \ \ \ \ \ \  \ \ \ \ \ \ \ \ \ (5)
\end{eqnarray*}
Then by (3), we can select $\delta$ small enough such that for each $\varepsilon\leq \delta,$ we have
$$\left\|\sup_{t\leq s\leq t+\varepsilon\wedge\tau}|Y_s^{t+\varepsilon\wedge\tau}|\right\|_\infty\leq \frac{1}{2\gamma}.\eqno(6)$$
Then by (5), (6) and (2), for each $\varepsilon\leq \delta,$ we have
\begin{eqnarray*}
\frac{1}{2}E\left[\int_{t}^{t+\varepsilon\wedge\tau}|{Z}_r^{t+\varepsilon\wedge\tau}|^2dr|{\cal{F}}_t\right]
&\leq&2E\left[\int_{t}^{t+\varepsilon\wedge\tau}\alpha_r|{Y}_r^{t+\varepsilon\wedge\tau}|dr|{\cal{F}}_t\right]
+2\beta E\left[\int_{t}^{t+\varepsilon\wedge\tau}|{Y}_r^{t+\varepsilon\wedge\tau}|^2dr|{\cal{F}}_t\right]\\
&\leq&2e^{\beta\varepsilon}\left\|\int_t^{t+\varepsilon}\alpha_rdr\right\|_\infty^2
+2\beta \varepsilon e^{2\beta\varepsilon}\left\|\int_t^{t+\varepsilon}\alpha_rdr\right\|_\infty^2\\
&=&2e^{\beta\varepsilon}\left\|\left(\int_t^{t+\varepsilon}\alpha_rdr\right)^2\right\|_\infty
+2\beta \varepsilon e^{2\beta\varepsilon}\left\|\left(\int_t^{t+\varepsilon}\alpha_rdr\right)^2\right\|_\infty\\
&\leq&2e^{\beta\varepsilon}\varepsilon\left\|\int_t^{t+\varepsilon}\alpha_r^2dr\right\|_\infty
+2\beta 2e^{2\beta\varepsilon}\varepsilon^2\left\|\int_t^{t+\varepsilon}\alpha_r^2dr\right\|_\infty.
\end{eqnarray*}
From above inequality and Remark 2.1, the proof is complete. $\Box$\\

Inspired by  Lepeltier and San Martin (1997, Lemma 1) and Fan et al. (2011, Lemma 3), we have the following Lemma 3.2.\\

\textbf{Lemma 3.2} Let $f(\cdot):{\textbf{R}}^k\longmapsto {\textbf{R}}$ with $k\in \textbf{N}$ is a continuous function with convex growth, that is, there exist constant $K\geq0, a\geq0$ and a strictly increasing convex function $\varphi:\ {\textbf{R}}^+\longrightarrow {\textbf{R}}^+$ with $\varphi(0)=\varphi(0+)=0,$ such that
$$|f(x)|\leq a+K\varphi(|x|).\eqno(7)$$
Then the following sequence of functions $f_n:$
$${f}_n(x)=\inf\{f(u)+\frac{n}{2}K\varphi(2|u-x|):{u\in{\textbf{Q}}^k}\}.\eqno(8)$$
is well-defined, for $n\geq1$ and we have

(i) $|{f}_n(x)|\leq a+\frac{1}{2}K\varphi(2|x|);$

(ii) ${f}_n(x)\nearrow$ as $n\longrightarrow\infty;$

(iii) ${f}_n(x)\longrightarrow f(x)$ as $n\longrightarrow\infty.$\\

\emph{Proof.} Clearly, we only need prove the case $K>0.$ The method of proof derives from Lepeltier and San Martin (1997, Lemma 1) and Fan et al. (2011, Lemma 3). By (7), (8), the convexity of $\varphi$ and $\varphi(0)=0,$ we have
\begin{eqnarray*}
f_n(x)&\leq& f(x)\leq a+K \varphi(|x|)\leq a+\frac{1}{2}K\varphi(2|x|)\\
f_n(x)&\geq& \inf\{f(u)+\frac{1}{2}K\varphi(2|u-x|):{u\in {\textbf{Q}}^k}\}\\
&\geq& \inf\{-a-K\varphi(|u-x+x|)+\frac{1}{2}K\varphi(2|u-x|):{u\in {\textbf{Q}}^k}\}\\
&\geq& -a-\frac{1}{2}K\varphi(2|x|)
\end{eqnarray*}
Thus, we get (i) and for $n\geq1,$ $f_n$ is well-defined. (ii) can be obtained from (8), directly. Let $\{x_n\}_{n\geq1}$ be a sequence such that $x_n\rightarrow x$ as $n\rightarrow \infty.$ By (7), (8) and the convexity of $\varphi$, for $n\geq1,$ we can take $u_n$ such that
\begin{eqnarray*}
\ \ \ \ \ \ \ \ \ \ \  \ \ \ \ f_n(x_n)&\geq& f(u_n)+\frac{n}{2}K\varphi(2|u_n-x_n|)-\frac{1}{n}\\
&\geq& -a-K\varphi(|u_n-x_n+x_n|)+\frac{n}{2}K\varphi(2|u_n-x_n|)-\frac{1}{n}\\
&\geq& -a-\frac{1}{2}K\varphi(2|x_n|)+\frac{n-1}{2}K\varphi(2|u_n-x_n|)-\frac{1}{n}
\ \ \ \ \ \ \ \ \ \ \ \ \ \ \  \ \ \ \ \ \ \ \ \ \ \ \ \ (9)
\end{eqnarray*}
Then by (i), we have
\begin{eqnarray*}
\frac{n-1}{2}K\varphi(2|u_n-x_n|)\leq2a+K\varphi(2|x_n|)+\frac{1}{n}
\end{eqnarray*}
Then we have $\limsup_{n\rightarrow\infty}\frac{n-1}{2}K\varphi(2|u_n-x_n|)\leq\infty.$ Since $\varphi$ is continuous and strictly increasing with $\varphi(0)=0,$ we have $u_n\rightarrow x$ as $n\rightarrow \infty.$
Then by (8), (9) and continuity of $f$, we have
\begin{eqnarray*}
\limsup_{n\rightarrow\infty}f_n(x_n)&\leq& \limsup_{n\rightarrow\infty}f(x_n)=f(x);\\
\liminf_{n\rightarrow\infty}f_n(x_n)&\geq& \liminf_{n\rightarrow\infty}f(u_n)=f(x).
\end{eqnarray*}
The proof is complete. $\Box$\\

\textbf{Remark 3.1} In Lemma 3.2, the cases $\varphi(x)=|x|$ and $\varphi(x)=|x|^b$ for $b\geq1,$ have been proved by Lepeltier and San Martin (1997, Lemma 1) and Fan et al. (2011, Lemma 3), respectively.\\

Inspired by Fan and Jiang (2010, Proposition 3) and Fan et al. (2011, Proposition 2), we can get the following Lemma 3.3 from Lemma 3.2. \\

\textbf{Lemma 3.3} Let $g$ satisfies assumption (A), for any $(y,x,q)\in{\mathbf{R}}\times {\mathbf{R}}^{\mathit{m}}\times {\mathbf{R}}^{\mathit{m}}$, there exists a non-negative process sequence $\{(\psi^n_t)_{t\in[0,T]}\}_{n=1}^{\infty}$ depending on $(y,x,q),$ satisfying for each $t\in[0,T],$ $\lim\limits_{n\rightarrow\infty}\psi^n_t=0$ and for each $n\geq1,$
$$|\psi^n_t|\leq4\alpha_t+M,\eqno(10)$$
for some constant $M>0$ only depending on $y,x,q,\gamma,\nu,$ such that for each $n\geq1$ and $(t,\bar{y},\bar{z},\bar{x})\in [0,T]\times{\mathbf{R}}\times {\mathbf{R}}^{d+m}$, we have
$$|g(t,\bar{y},\bar{z}+\sigma^\ast(t,\bar{x})q)-g(t,y,\sigma^\ast(t,x)q)|\leq \frac{n}{2}\varphi(2|\bar{y}-y|)+2 n\lambda(|\bar{z}|^2+|\bar{x}-x|^2)+{\psi}^n_t,$$
where $\lambda:=\gamma+2\gamma|q|^2|\nu|^2$ and $\nu$ is the constant in (H2).\\

\emph{Proof.} For any $(y,x,q)\in{\mathbf{R}}\times {\mathbf{R}}^{\mathit{m}}\times {\mathbf{R}}^{\mathit{m}}$, we set $$f(t,y,z,x):=g(t,y,z+\sigma^\ast(t,x)q).$$
Then by (A) and (H2), we have
\begin{eqnarray*}
 \ \ \ \ \ \ \ \  \ \ \ \ \ \ \ \ |f(t,y,z,x)|&\leq& \alpha_t+\varphi(|y|)+\frac{\gamma}{2}|z+\sigma^\ast(t,x)q|^2\\
&\leq&\alpha_t+\varphi(|y|)+\gamma|z|^2+2\gamma|q|^2|\nu|^2(1+|x|^2)\\
&\leq&\hat{\alpha}_t+\varphi(|y|)+\lambda(|z|^2+|x|^2).\ \ \ \ \ \ \ \ \ \ \ \ \ \ \ \  \ \ \ \ \ \ \ \ \ \ \ \  \ \ \ \ \ \ \ \  \ \ \ \ \ \ \ \ (11)
\end{eqnarray*}
where $\hat{\alpha}_t=\alpha_t+2\gamma|q|^2|\nu|^2$ and $\lambda=\gamma+2\gamma|q|^2|\nu|^2.$ We set
$${f}_n^1(t,y,z,x):=\sup\{f(t,u,v,w)-(\frac{n}{2}\varphi(2|u-y|)+2 n\lambda(|v-z|^2+|w-x|^2)):{(u,v,w)\in \textbf{Q}^{1+d+n}}\}.$$
$${f}_n^2(t,y,z,x):=\inf\{f(t,u,v,w)+(\frac{n}{2}\varphi(2|u-y|)+2 n\lambda(|v-z|^2+|w-x|^2)):{(u,v,w)\in \textbf{Q}^{1+d+n}}\}.$$
and
$${\psi}^1_n(t):={f}_n^1(t,y,0,x)\ \  \textrm{and}\ \ {\psi}^2_n(t):={f}_n^2(t,y,0,x).$$
By Lemma 3.2, we can deduce
$$|{\psi}^1_n(t)|\leq\hat{\alpha}_t+\frac{1}{2}\varphi(2|y|)+2\lambda|x|^2\ \  \textrm{and}\ \
|{\psi}^2_n(t)|\leq\hat{\alpha}_t+\frac{1}{2}\varphi(2|y|)+2\lambda|x|^2,\eqno(12)$$
and $$\lim_{n\rightarrow\infty}{\psi}^1_n(t)=\lim_{n\rightarrow\infty}{\psi}^2_n(t)=f(t,y,0,x).\eqno(13)$$
By the definition of $f^1_n$ and $f^2_n$, we also have
$$f(t,\bar{y},\bar{z},\bar{x})-f(t,y,0,x)\leq (\frac{n}{2}\varphi(2|\bar{y}-y|)+2 n\lambda(|\bar{z}|^2+|\bar{x}-x|^2))+{\psi}^1_n(t)-f(t,y,0,x);$$
$$f(t,\bar{y},\bar{z},\bar{x})-f(t,y,0,x)\geq -(\frac{n}{2}\varphi(2|\bar{y}-y|)+2 n\lambda(|\bar{z}|^2+|\bar{x}-x|^2))+{\psi}^2_n(t)-f(t,y,0,x).$$
By setting
$${\psi}^n_t:=|{\psi}^1_n(t)-f(t,y,0,x)|+|{\psi}^2_n(t)-f(t,y,0,x)|,\eqno(14)$$
we have
$$|f(t,\bar{y},\bar{z},\bar{x})-f(t,y,0,x)|\leq \frac{n}{2}\varphi(|\bar{y}-y|)+2 n\lambda(|\bar{z}|^2+|\bar{x}-x|^2)+{\psi}^n_t.\eqno(15)$$
By (11)-(15), we can complete this proof.  $\Box$\\

The following Theorem 3.1 is the main result of this paper.\\

\textbf{Theorem 3.1} Let $g$ satisfies assumption (A), for each $(y,x,q)\in{\mathbf{R}}\times{\mathbf{R}}^{\mathit{m}}\times {\mathbf{R}}^{\mathit{m}}$ and almost every $t\in[0,T[$, we have

(i) if $\{\alpha_t\}_{t\in[0,T]}$ satisfies assumption (B) and $\|\alpha_t\|_{{\cal{S}}^{{1}}}<\infty,$ then
$$\ g\left(t,y,\sigma^\ast(t,x)q\right)+q\cdot b(t,x)
=L^1-\lim\limits_{\varepsilon\rightarrow0^+}\frac{1}{\varepsilon}
\left(Y_t^{t+\varepsilon\wedge\tau}-y\right);$$

(ii) if $\|\alpha_t\|_\infty<\infty,$ then
$$P-a.s.,\ \ \ g\left(t,y,\sigma^\ast(t,x)q\right)+q\cdot b(t,x)
=\lim\limits_{\varepsilon\rightarrow0^+}\frac{1}{\varepsilon}
\left(Y_t^{t+\varepsilon\wedge\tau}-y\right).$$
where $\tau=\inf\{s\geq0:|X_{t+s}^{t,x}|> C_0\}\wedge (T-t)$ for constant $C_0>|x|$ and $(Y_s^{t+\varepsilon\wedge\tau},Z_s^{t+\varepsilon\wedge\tau})$ is an arbitrary solution of BSDE with parameter $(g,t+\varepsilon\wedge\tau,y+q\cdot(X_{t+\varepsilon\wedge\tau}^{t,x}-x)).$\\

\textit{Proof.} Let $(y,x,q)\in{\mathbf{R}}\times{\mathbf{R}}^{\mathit{m}}\times {\mathbf{R}}^{\mathit{m}}.$  For a constant $C_0>|x|,$ we define the following stopping time:
$$\tau:=\inf\left\{s\geq0:|X_{t+s}^{t,x}|> C_0\right\}\wedge (T-t).$$
By the continuity of $X_{t+s}^{t,x}$, we have
$0<\tau\leq T-t$ and for $s\in[0,{t+\varepsilon\wedge\tau}],$
$$|X_{s}^{t,x}|\leq C_0.\eqno(16)$$

For $\varepsilon\in]0,T-t]$, let $\left(Y_s^{t+\varepsilon\wedge\tau},Z_s^{t+\varepsilon\wedge\tau}\right)$ be a solution of BSDE with parameter $(g,t+\varepsilon\wedge\tau,y+q\cdot(X_{t+\varepsilon\wedge\tau}^{t,x}-x))$ and set for $s\in[t,{t+\varepsilon\wedge\tau}],$
$$
\tilde{Y}_s^{{t+\varepsilon\wedge\tau}}:=Y_s^{{t+\varepsilon\wedge\tau}}-(y+q\cdot(X_{s}^{t,x}-x)),\ \ \ \
\tilde{Z}_s^{t+\varepsilon\wedge\tau}:=Z_s^{t+\varepsilon\wedge\tau}-\sigma^\ast(t,X_{s}^{t,x})q.\eqno(17)
$$
Applying It\^{o} formula to $\tilde{Y}_s^{t+\varepsilon\wedge\tau}$ for $s\in[t,{t+\varepsilon\wedge\tau}],$ we have
\begin{eqnarray*}
\tilde{Y}_s^{t+\varepsilon\wedge\tau}&=&\int_{s}^{{t+\varepsilon\wedge\tau}}\left(g(r,\tilde{Y}_r^{t+\varepsilon\wedge\tau}+y
+q\cdot(X_{r}^{t,x}-x),
\tilde{Z}_r^{t+\varepsilon\wedge\tau}+\sigma^\ast(r,X_{r}^{t,x})q)+q\cdot b(r,X_{r}^{t,x})\right)dr\\
&&-\int_{s}^{t+\varepsilon\wedge\tau}\tilde{Z}_r^{t+\varepsilon\wedge\tau}dB_r.\ \ \ \ \ \ \ \ \ \ \ \ \ \ \ \ \ \ \ \ \ \ \ \ \ \ \ \ \ \ \ \ \ \ \ \ \ \ \ \ \ \ \ \ \ \ \ \ \ \ \ \ \ \ \ \ \ \ \ \ \ \ \ \ \ \ \ \ \ \ \ \ \ \ \ (18)
\end{eqnarray*}
Set $\tilde{g}(r,\tilde{y},\tilde{z}):=g(r,\tilde{y}+y
+q\cdot(X_{r}^{t,x}-x),
\tilde{z}+\sigma^\ast(r,X_{r}^{t,x})q)+q\cdot b(r,X_{r}^{t,x}).$ Clearly, $\tilde{g}(r,\tilde{y},\tilde{z})$ satisfies (i) and (ii) in (A). Moreover, by (A), the convexity of $\varphi$, (H1) and (16), we have for $r\in[0,{t+\varepsilon\wedge\tau}],$
\begin{eqnarray*}
|\tilde{g}(r,\tilde{y},\tilde{z})|&=&|g(r,\tilde{y}+y+q\cdot(X_{r}^{t,x}-x),\tilde{z}+\sigma^\ast(r,X_{r}^{t,x})q)+q\cdot b(r,X_{r}^{t,x})|\\
&\leq &\alpha_t+\varphi(|\tilde{y}+y+q\cdot(X_{r}^{t,x}-x)|)+\frac{\gamma}{2}|\tilde{z}+\sigma^\ast(r,X_{r}^{t,x})q|^2+|q\cdot b(r,X_{r}^{t,x})|
\\
&\leq&\alpha_t+\frac{1}{2}\varphi(2|\tilde{y}|)+\frac{1}{2}\varphi(2|y+q\cdot(X_{r}^{t,x}-x)|)+\gamma|\tilde{z}|^2
+\gamma|\sigma^\ast(r,X_{r}^{t,x})q|^2+|q\cdot b(r,X_{r}^{t,x})|\\
&\leq&\tilde{\alpha}_t+\frac{1}{2}\varphi(2|\tilde{y}|)+\gamma|\tilde{z}|^2
\end{eqnarray*}
where $$\tilde{\alpha}_t=\alpha_t+\frac{1}{2}\varphi(2(|y|+|q|C_0+|q||x|))+\gamma\nu^2q^2(1+C_0)^2+q\nu(1+C_0).\eqno(19)$$ Then, we get that $\tilde{g}$ also satisfy (iii) in (A). By (18), we get $(\tilde{Y}_s^{t+\varepsilon\wedge\tau},\tilde{Z}_s^{t+\varepsilon\wedge\tau})$ is a solution of BSDE with parameter $(\tilde{g},t+\varepsilon\wedge\tau,0)$ in $[t,t+\varepsilon\wedge\tau].$ If $\{\alpha_t\}_{t\in[0,T]}$ satisfies (B), then by (19) and (3),  we have,
$$\lim\limits_{\varepsilon\rightarrow0^+}\left\|\sup_{t\leq s\leq t+\varepsilon\wedge\tau}|\tilde{Y}_{s}^{t+\varepsilon\wedge\tau}|\right\|_\infty=0.
\eqno(20)$$
Moreover, if $\|\alpha_t\|_\infty<\infty,$ by (19) and (4),  we have
$$\sup_{t\leq s\leq t+\varepsilon\wedge\tau}|\tilde{Y}_{s}^{t+\varepsilon\wedge\tau}|\leq e^{\beta \varepsilon}\varepsilon\|\tilde{\alpha}_t\|_\infty. \eqno(21)$$
Set
\begin{eqnarray*}\
M^{\varepsilon,\tau}_t&:=&\frac{1}{\varepsilon}E\left[\int_{t}^{t+\varepsilon\wedge\tau}g(r,\tilde{Y}_r^{t+\varepsilon\wedge\tau}
+y+q\cdot(X_{r}^{t,x}-x),
\tilde{Z}_r^{t+\varepsilon\wedge\tau}+\sigma^\ast(r,X_{r}^{t,x})q)dr|{\cal{F}}_t\right]\\
P^{\varepsilon,\tau}_t&:=&\frac{1}{\varepsilon}E\left[\int_{t}^{t+\varepsilon\wedge\tau}g(r,y,
\sigma^\ast(r,x)q)dr|{\cal{F}}_t\right],\\
U^{\varepsilon,\tau}_t&:=&\frac{1}{\varepsilon}E\left[\int_{t}^{t+\varepsilon\wedge\tau}q\cdot b(r,X_{r}^{t,x})dr|{\cal{F}}_t\right],
\end{eqnarray*}
Then by (17) and (18), we have
\begin{eqnarray*}
&&\frac{1}{\varepsilon}\left(Y_t^{t+\varepsilon\wedge\tau}-y
\right)
-g(t,y,\sigma^\ast(t,x)q)-q\cdot b(t,x)\\
&=&\frac{1}{\varepsilon}\tilde{Y}_t^{t+\varepsilon\wedge\tau}
-g(t,y,\sigma^\ast(t,x)q)-q\cdot b(t,x)\\
&=&\left(M^{\varepsilon,\tau}_t-P^{\varepsilon,\tau}_t\right)
+\left(P^{\varepsilon,\tau}_t-g(t,y,\sigma^\ast(t,x)q)\right)+\left(U^{\varepsilon,\tau}_t-q\cdot b(t,x)\right).\ \ \ \  \ \ \ \  \ \  \ \ \ \  \ \  \ \ \ \  \ \  \ \ \ \ \ \ \ \ (22)
\end{eqnarray*}
By Jensen inequality, we have
\begin{eqnarray*}
 &&|M^{\varepsilon,\tau}_t-P^{\varepsilon,\tau}_t|\\
&\leq&\frac{1}{\varepsilon}E\left[\left|\int_{t}^{t+\varepsilon\wedge\tau}\left(g(r,\tilde{Y}_r^{t+\varepsilon\wedge\tau}+y
+q\cdot(X_{r}^{t,x}-x),
\tilde{Z}_r^{t+\varepsilon\wedge\tau}+\sigma^\ast(r,X_{r}^{t,x})q)-g(r,y,
\sigma^\ast(r,x)q)\right)dr\right||{\cal{F}}_t\right]\\
&\leq&\frac{1}{\varepsilon}E\left[\int_{t}^{t+\varepsilon\wedge\tau}\left|g(r,\tilde{Y}_r^{t+\varepsilon\wedge\tau}+y
+q\cdot(X_{r}^{t,x}-x),
\tilde{Z}_r^{t+\varepsilon\wedge\tau}+\sigma^\ast(r,X_{r}^{t,x})q)-g(r,y,
\sigma^\ast(r,x)q)\right|dr|{\cal{F}}_t\right]\\
\end{eqnarray*}
Then by Lemma 3.3, there exists a non-negative process sequence $\{(\psi^n_t)_{t\in[0,T]}\}_{n=1}^{\infty}$ depending on $(y,x,q),$ such that for each $n\in \textbf{N},$
\begin{eqnarray*}
&&|M^{\varepsilon,\tau}_t-P^{\varepsilon,\tau}_t|\\
&\leq&\frac{1}{\varepsilon}E\left[\int_{t}^{t+\varepsilon\wedge\tau}\left|\frac{n}{2}\varphi(2|\tilde{Y}_r^{t+\varepsilon\wedge\tau}
+q(X_{r}^{t,x}-x)|)+2n\lambda(|\tilde{Z}_r^{t+\varepsilon\wedge\tau}|^2+|X_{r}^{t,x}-x|^2)
+\psi_r^n\right|dr|{\cal{F}}_t\right]\\
&\leq&\frac{n}{4\varepsilon}E\left[\int_{t}^{t+\varepsilon\wedge\tau}\varphi(4|\tilde{Y}_r^{t+\varepsilon\wedge\tau}|)dr|{\cal{F}}_t\right]
+\frac{2}{\varepsilon}n\lambda E\left[\int_{t}^{t+\varepsilon\wedge\tau}|\tilde{Z}_r^{t+\varepsilon\wedge\tau}|^2dr|{\cal{F}}_t\right]\\
&&+\frac{n}{4\varepsilon}E\left[\int_{t}^{t+\varepsilon\wedge\tau}\varphi(4|q||X_{r}^{t,x}-x|)dr|{\cal{F}}_t\right]
+\frac{2}{\varepsilon}n\lambda E\left[\int_{t}^{t+\varepsilon\wedge\tau}|X_{r}^{t,x}-x|^2dr|{\cal{F}}_t\right]\\
&&+\frac{1}{\varepsilon}E\left[\int_{t}^{t+\varepsilon}|\psi_r^n|dr|{\cal{F}}_t\right]\ \ \ \ \ \ \ \ \ \ \ \ \ \ \ \ \ \ \ \ \ \ \ \ \ \ \ \ \ \ \ \ \ \ \ \ \ \ \ \  \ \ \ \ \ \ \ \ \ \ \ \ \ \ \ \ \ \ \ \ \ \ \ \ \ \ \ \ \ \ \  \ \ \ \ \ \ \ \ \ \ (23)
\end{eqnarray*}
where $\lambda:=\gamma+2\gamma|q|^2|\nu|^2$ and $\nu$ is the constant in (H2).\\

\textbf{Proof of (i):} If $\{\alpha_t\}_{t\in[0,T]}$ satisfies (B), then by (20), we have
$$\lim\limits_{\varepsilon\rightarrow0^+}\frac{1}{\varepsilon}E\left[\int_{t}^{t+\varepsilon\wedge\tau}
\varphi(4|\tilde{Y}_r^{t+\varepsilon\wedge\tau}|)dr\right]\leq\lim\limits_{\varepsilon\rightarrow0^+}
\varphi(4\|\sup_{t\leq s\leq t+\varepsilon\wedge\tau}|\tilde{Y}_{s}^{t+\varepsilon}|\|_\infty)=0,\eqno(24)$$
and
$$\lim\limits_{\varepsilon\rightarrow0^+}\frac{1}{\varepsilon}E\left[\int_{t}^{t+\varepsilon\wedge\tau}
\varphi(4|\tilde{Y}_r^{t+\varepsilon\wedge\tau}|)dr|{\cal{F}}_t\right]\leq\lim\limits_{\varepsilon\rightarrow0^+}
\varphi(4\|\sup_{t\leq s\leq t+\varepsilon\wedge\tau}|\tilde{Y}_{s}^{t+\varepsilon}|\|_\infty)=0,\eqno(25)$$
By (19) and Lemma 3.1, we have
$$\lim\limits_{\varepsilon\rightarrow0^+}\frac{1}{\varepsilon} E\left[\int_{t}^{t+\varepsilon\wedge\tau}|\tilde{Z}_r^{t+\varepsilon\wedge\tau}|^2dr\right]=0\ \ \textrm{and}\ \ \lim\limits_{\varepsilon\rightarrow0^+}\frac{1}{\varepsilon} E\left[\int_{t}^{t+\varepsilon\wedge\tau}|\tilde{Z}_r^{t+\varepsilon\wedge\tau}|^2dr|{\cal{F}}_t\right]=0.\eqno(26)$$
By (16), Lebesgue dominated convergence theorem and the continuity of $X_{r}^{t,x}$ in $r$, we have,
$$\lim\limits_{\varepsilon\rightarrow0^+}\frac{1}{\varepsilon}E\left[\int_{t}^{t+\varepsilon\wedge\tau}\varphi(4|q||X_{r}^{t,x}-x|)dr|{\cal{F}}_t\right]
=E\left[\lim\limits_{\varepsilon\rightarrow0^+}\frac{1}{\varepsilon}\int_{t}^{t+\varepsilon\wedge\tau}
\varphi(4|q||X_{r}^{t,x}-x|)dr|{\cal{F}}_t\right]=0,\eqno(27)$$
and
$$\lim\limits_{\varepsilon\rightarrow0^+}\frac{1}{\varepsilon}E\left[\int_{t}^{t+\varepsilon\wedge\tau}
\varphi(4|q||X_{r}^{t,x}-x|)dr\right]
=E\left[\lim\limits_{\varepsilon\rightarrow0^+}\frac{1}{\varepsilon}\int_{t}^{t+\varepsilon\wedge\tau}
\varphi(4|q||X_{r}^{t,x}-x|)dr\right]
=0.\eqno(28)$$
Since $\|\alpha_t\|_{{\cal{S}}^{{1}}}<\infty,$ by (10) and Lebesgue dominated convergence theorem, we have $$\lim_{n\rightarrow\infty}\|\psi_t^n\|_{{\cal{H}}^1}=E\int_0^T\lim_{n\rightarrow\infty}\psi_t^ndt=0.\eqno(29)$$
Thus by Fatou Lemma and Fubini Theorem, we have
$$\int_{0}^T\liminf_{n\rightarrow\infty}E|\psi_r^n|dr\leq \liminf_{n\rightarrow\infty}\int_{0}^TE|\psi_r^n|dr=
\liminf_{n\rightarrow\infty}E\int_{0}^T|\psi_r^n|dr
=\lim\limits_{n\rightarrow\infty}
\|\psi^n_t\|_{\cal{H}^{\mathrm{1}}}=0.$$
Then, for almost every $t\in[0,T[$, we have
$$\liminf_{n\rightarrow\infty}E|\psi_t^n|=0.\eqno(30)$$
Taking expectation on both sides of (23), then by (24), (26), (28), Fubini Theorem, Lemma 2.1 and (30), we get that for almost every $t\in[0,T[$,
\begin{eqnarray*}
\ \ \ \ \ \lim\limits_{\varepsilon\rightarrow0^+}E|M^{\varepsilon,\tau}_t
-P^{\varepsilon,\tau}_t|\leq\liminf_{n\rightarrow\infty}\lim\limits_{\varepsilon\rightarrow0^+}\frac{1}{\varepsilon}
E\left[\int_{t}^{t+\varepsilon}|\psi_r^n|dr\right]
&=&\liminf_{n\rightarrow\infty}\lim\limits_{\varepsilon\rightarrow0^+}\frac{1}{\varepsilon}
\int_{t}^{t+\varepsilon}E|\psi_r^n|dr\\
&=&\liminf_{n\rightarrow\infty}E|\psi_t^n|\\
&=&0.\ \ \ \ \ \ \ \ \ \ \ \ \ \  \ \ \ \ \ \ \ \ \ \  \ \ \ \ \ \ \ \ \ \ \ (31)
\end{eqnarray*}
Since $\|\alpha_t\|_{{\cal{S}}^{{1}}}<\infty,$ then by the same argument as (20) and (19) in Zheng (2014b), we have for almost every $t\in[0,T[$,
$$\lim\limits_{\varepsilon\rightarrow0^+}E|P^{\varepsilon,\tau}_t-g(t,y,\sigma^\ast(t,x)q)|=0,\ \
\textrm{and} \ \ \lim\limits_{\varepsilon\rightarrow0^+}E|U^{\varepsilon,\tau}_t-q\cdot b(t,x)|=0,\eqno(32)$$
respectively. Then, by (22) , (31) and (32), we get (i).

\textbf{Proof of (ii): } If $\|\alpha_t\|_\infty<\infty,$ then $\{\alpha_t\}_{t\in[0,T]}$ satisfies (B).
By (23), (25)-(27), Zheng (2014a, Lemma 3.3) and Lemma 3.3, we get that for almost every $t\in[0,T[$,
$$ P-a.s.,\ \ \lim\limits_{\varepsilon\rightarrow0^+}|M^{\varepsilon,\tau}_t
-P^{\varepsilon,\tau}_t|\leq\liminf_{n\rightarrow\infty}\lim\limits_{\varepsilon\rightarrow0^+}\frac{1}{\varepsilon}
E\left[\int_{t}^{t+\varepsilon}|\psi_r^n|dr|{\cal{F}}_t\right]
=\lim_{n\rightarrow\infty}|\psi_t^n|
=0.\eqno(33)$$
by (A), (H2) and Zheng (2014a, Lemma 3.3), we have for almost every $t\in[0,T[$,
$$P-a.s.,\ \ \lim\limits_{\varepsilon\rightarrow0^+}|P^{\varepsilon,\tau}_t-g(t,y,\sigma^\ast(t,x)q)|=0.\eqno(34)$$
By (20) in Zheng (2014a), we have, for almost every $t\in[0,T[$,
$$P-a.s.,\ \ \lim\limits_{\varepsilon\rightarrow0^+}|U^{\varepsilon,\tau}_t-q\cdot b(t,x)|=0.\eqno(35)$$
By (22), (33)-(35), we get (ii). The proof is completed.   $\Box$\\

If $\|\alpha_t\|_\infty<\infty,$ then by (17), (19), (21), Theorem 3.1 and the same arguments as Zheng (2014a, Theorem 3.2 and Theorem 3.3), we have the following Corollary 3.1. We omit its proof.\\

\textbf{Corollary 3.1} Let $g$ satisfy assumption (A) and $\|\alpha_t\|_\infty<\infty,$  then for each $(y,x,q)\in{\mathbf{R}}\times{\mathbf{R}}^{\mathit{n}}\times {\mathbf{R}}^{\mathit{n}}$ and $p>0$, we have
$$g\left(t,y,\sigma^\ast(t,x)q\right)+q\cdot b(t,x)
=L^p -\lim\limits_{\varepsilon\rightarrow0^+}\frac{1}{\varepsilon}
\left(Y_t^{t+\varepsilon\wedge\tau_t}-y\right),\ \ \ a.e.\ t\in[0,T[,$$
and
$$g\left(t,y,\sigma^\ast(t,x)q\right)+q\cdot b(t,x)
={\cal{H}}^p-\lim\limits_{\varepsilon\rightarrow0^+}\frac{1}{\varepsilon}
\left(Y_t^{t+\varepsilon\wedge\tau_t}-y\right),$$
where $\tau_t=\inf\{s\geq0:|X_{t+s}^{t,x}|> C_t\}\wedge (T-t)$ for constant $C_t>|x|$ and $(Y_s^{t+\varepsilon\wedge\tau_t},Z_s^{t+\varepsilon\wedge\tau_t})$ is an arbitrary solution of BSDE with parameter $(g,t+\varepsilon\wedge\tau_t,y+q\cdot(X_{t+\varepsilon\wedge\tau_t}^{t,x}-x)).$\\

Let $q=z, b(t,x)=0, \sigma(t,x)=1, x=0$ in Theorem 3.1, then we have the following Corollary 3.3, immediately.\\

\textbf{Corollary 3.2} Let $g$ satisfies assumption (A), then for each $(y,z)\in{\mathbf{R}}\times {\mathbf{R}}^{\mathit{d}}$ and almost every $t\in[0,T[$, we have

(i) if $\{\alpha_t\}_{t\in[0,T]}$ satisfies assumption (B) and $\|\alpha_t\|_{{\cal{S}}^{{1}}}<\infty,$ then
$$g\left(t,y,z\right)
=L^1-\lim\limits_{\varepsilon\rightarrow0^+}\frac{1}{\varepsilon}
\left(Y_t^{t+\varepsilon\wedge\tau}-y\right);$$

(ii) if $\|\alpha_t\|_\infty<\infty,$ then
$$P-a.s.,\ \ \ g\left(t,y,z\right)
=\lim\limits_{\varepsilon\rightarrow0^+}\frac{1}{\varepsilon}
\left(Y_t^{t+\varepsilon\wedge\tau}-y\right),$$
where $\tau=\inf\{s\geq0:|B_{t+s}-B_t|> C_0\}\wedge (T-t)$ for constant $C_0>|x|$ and $(Y_s^{t+\varepsilon\wedge\tau},Z_s^{t+\varepsilon\wedge\tau})$ is an arbitrary solution of BSDE with parameter $(g,t+\varepsilon\wedge\tau,y+z\cdot(B_{t+\varepsilon\wedge\tau}-B_t)).$\\

By Corollary 3.2 and a simple discussion, we can get the following converse comparison theorem. We omit its proof here.\\

{\textbf{Corollary 3.3}} Let generators $g_1$ and $g_2$ both satisfy assumption (A) in which $\{\alpha_t\}_{t\in[0,T]}$ satisfies assumption (B) and for any stopping time $\tau \in]0,T],\ \xi\in L^\infty({\mathcal{F}}_\tau),$ BSDEs with parameter $(g_1,\tau,\xi)$ and $(g_2,\tau,\xi)$ exist solutions $(Y_{t\wedge\tau}^{\tau,1},Z_{t\wedge\tau}^{\tau,1})$ and $(Y_{t\wedge\tau}^{\tau,2},Z_{t\wedge\tau}^{\tau,2})$, respectively, such that $\forall t\in[0,T],$
$$P-a.s.,\ \ Y_{t\wedge\tau}^{\tau,1}\geq Y_{t\wedge\tau}^{\tau,2},\eqno(36)$$
then for each $(y,z)\in{\mathbf{R}}\times{\mathbf{R}}^{\mathit{d}},$ and almost every $t\in[0,T[$, we have
$$P-a.s.,\  \ g_1(t,y,z)\geq g_2(t,y,z).$$

\textbf{Remark 3.2} Fan et al. (2011, Theorem 1) establishes a representation theorem for generator of BSDE in ${\cal{H}}^p$ for $1\leq p<2,$ when generator is monotonic and polynomial growth in $y$ ($\varphi(\cdot)=|\cdot|^b$) and linear growth in $z$ and $\alpha_t\in {\cal{H}}^{2b}_T({\textbf{R}}^d), b\geq1.$ One can see (i) in Theorem 3.1 is a representation theorem in $L^1$. In fact, under the same conditions, using the method of proof of Theorem 3.1 and some results in the Fan $\&$ Jiang (2010), we can establish a representation theorem in ${\cal{H}}^1$.\\

\textbf{Remark 3.3} In fact, in the spirit of Briand et al. (2000), Jiang (2005c, 2006, 2008) and Jia (2008). we can use the representation theorem for generator obtained in this paper to study the properties of such BSDEs. Moreover, when solution of such BSDE is unique and generator $g(\cdot,\cdot,0)\equiv0$, we also can apply this representation theorem to study the properties of $g$-expectation induced by such BSDEs.


\begin{thebibliography}{99}
\addtolength{\itemsep}{-0.8 em}
\bibitem[1]{1} Briand, P., Coquet, F., Hu, Y., M\'emin, J., Peng, S., 2000. A converse
comparison theorem for BSDEs and related properties of
$g$-expectation. Electron. Comm. Probab. 5, 101-117.
\bibitem[2]{2} Briand, P, Hu, Y. 2008. Quadratic BSDEs with convex generators and unbounded terminal conditions. Probability Theory and Related Fields, , 141(3-4): 543-567.
\bibitem[3]{3} Fan, S., Jiang, L., 2010. A representation theorem for generators of BSDEs with continuous linear-growth
generators in the space of processes. Journal of Computational and Applied mathematics, 235(3),
686-695.
\bibitem[4]{4} Fan, S., Jiang, L., Xu, Y., 2011. Representation theorem for generators of BSDEs with monotonic and polynomial-growth generators in the space of processes. Electronic Journal of Probability, 16(27), 830-834.
\bibitem[5]{5} Hewitt, E., Stromberg, K.R., 1978. Real and Abstract Analysis. Springer-Verlag, New York.
\bibitem[6]{6} Jia, G., 2008. Backward stochastic differential equations, $g$-expectations and related
semilinear PDEs. PH.D Thesis, ShanDong University, China, 2008.
\bibitem[7]{7}Jiang, L., 2005a. Representation theorems for generators of backward stochastic differential equations. Comptes Rendus Mathematique, 340(2), 161-166.
\bibitem[8]{8} Jiang, L., 2005b. Representation theorems for
generators of backward stochastic differential equations and their
applications. Stochastic Process. Appl. 115 (12), 1883-1903.
\bibitem[9]{9} Jiang, L., 2005c. Nonlinear expectation¡ª$g$-expectation theory and its
applications in finance. Ph.D Thesis. ShanDong University, China.
\bibitem[10]{10} Jiang, L. 2006. Limit theorem and uniqueness theorem of backward stochastic differential equations. Science in China Series A: Mathematics, 49(10), 1353-1362.
\bibitem[11]{11} Jiang, L., 2008. Convexity, translation invariance and subadditivity for
g-expectations and related risk measures. Annals of Applied
Probability 18 (1), 245-258.
\bibitem[12]{12} Liu, Y., Jiang, L., Xu, Y., 2008. A local limit theorem for solutions of
BSDEs with Mao' non-Lipschitz generator. Acta Math. Appli. Sinica,
English Series 24, 329-336.
\bibitem[13]{13} Pardoux, E., Peng, S., 1990. Adapted solution of backward stochastic
differential equations. Systems Control Letters 14, 51-61.
\bibitem[14]{14} Zheng, S., 2014a. Representation theorems for generators of quadratic BSDEs. arXiv:1405.4789.
\bibitem[15]{15} Zheng, S., 2014b. Representation theorems for generators of reflected BSDEs with continuous and linear-growth generators.  arXiv:1404.0226.


\end{thebibliography}
\end{document}